\documentclass[reqno,12pt]{amsart}
\usepackage[utf8]{inputenc}

\usepackage{graphics,amssymb,amsthm,amsmath}
\usepackage[english]{babel}
\usepackage{fullpage}
\usepackage{color}

\def\qpol{\overline{\mathcal M}_q({\mathfrak r})}

 \newtheorem{Thm}{Theorem}[section]
\newtheorem{Def}[Thm]{Definition} \newtheorem{Rem}[Thm]{Remark}
\newtheorem{Lem}[Thm]{Lemma} \newtheorem{Cor}[Thm]{Corollary}
\newtheorem{Prop}[Thm]{Proposition}

\newtheorem{Ex}[Thm]{Example}
\begin{document}
\fontsize{.496cm}{.496cm}\sf

\title[The center of a quantized nilpotent]{The center of ${\mathcal U}_q({\mathfrak n}_\omega)$. } 

\date{\today}
\author{Hans P. Jakobsen}
\address{
  Department of Mathematical Sciences\\ University of
Copenhagen\\Universitetsparken 5\\
   DK-2100, Copenhagen,
  Denmark} \email{jakobsen@math.ku.dk}
\begin{abstract}Based on \cite{jz}, we point to a new and very useful direction of approach to a general set of problems. We exemplify it here by obtaining the center of a localization of  ${\mathcal 
U}_q({\mathfrak 
n}_\omega)\subseteq {\mathcal U}^+_q({\mathfrak g})$  by the covariant elements 
(non-mutable elements). It is based on constructions and results from quantum cluster algebras. The 
non-zero complex parameter $q$ is mostly assumed not to be a root of unity, but our method also gives many details in case $q$ is a primitive root of unity.  Further, we use this to give a generalization to double Schubert Cell algebras and indicate that other families may also be handled. \end{abstract}
\maketitle

\section{Introduction}
The topics of quantum groups, quantized function algebras, quantized matrix 
algebras, and quantum cluster algebras have since long been seen to be 
intrinsically interwoven.

The groundbreaking research of  Drinfeld (\cite{drin1},\cite{drin2}) and Jimbo 
(\cite{jim1},\cite{jim2}) was followed by deep results of Lusztig 
(\cite{lu1},\cite{l})), Kashiwara (\cite{K1},\cite{K2}).  Then 
Levendorskii and Soibelman (\cite{lev2}, \cite{lev0}) and later  de Concini 
and Procesi (\cite{cp}) added the  quadratic algebra side to this distinguished 
family. With the advent of the cluster algebras of Fomin and Zelevinsky 
(\cite{fz}) 
and Berenstein-Zelevinsky quantized cluster algebras (\cite{bz}) many new 
dimensions were added to the function algebra side.

 Through many years, quantized function algebras have
attracted a lot of attention (\cite{cp}, \cite{cl},
\cite{dd}, \cite{frt}, \cite{ho-la}, \cite{jaz1}, \cite{jj} \cite{lev},
\cite{lev2}, \cite{lu1}, \cite{manin}, \cite{pw}, and many others). Many 
special examples were considered in the beginning, but also general families 
have more recently been considered (\cite{leclerc}). 

The current research has its focus on the quadratic algebra side. It utilizes 
fundamental results in (\cite{bz}) and (\cite{leclerc}). 

In certain families of examples (\cite{jz},\cite{jp}), it was seen that certain 
(signed) permutation matrices contained much information about the quantized 
matrix algebras. The topic of this article is to explain exactly the reason for 
that, while at the same time giving the full description of the centers. That 
we thereby also obtain an insight into  very algebraic properties of quantized 
function algebras, even specializing these to roots of unity, is clear, but 
will not be pursued in this article until the very last section. For now, and until then, we assume that $q$ is a 
not roots of unity. 

An important tool in our investigation is a family of quantized minors 
introduced by Berenstein and Zelevinsky (\cite{bz}).  Later C. Geiss, B. 
leclerc, and J. Schr\"oer (\cite{leclerc}) have modified these in a way that 
turns out to be exactly suitable for our needs.

\bigskip

Given an element $w$ in the Weyl group $W$ one may construct, using Lusztig 
(\cite{luz}), a family of elements $Z_1,\dots,Z_{\ell(w)}\subset {\mathcal 
U}_q^+$. It is a key result of Levendorskii and Soibelman's (\cite{lev2}, 
\cite{lev0}))  that, for 
$1\leq i<j\leq\ell(w)$, $$Z_{i}Z_{j}=q^{(\gamma_i,\gamma_j)}Z_{j}Z_{i} +
\textrm{terms involving only elements $Z_k$ with $i\leq k\leq
j$.}$$

These relations can be taken as the defining relations of ${\mathcal 
U}_q({\mathfrak n}_\omega)$. 

Procesi and de Concini later reproved this result and introduced the associated 
quasi-polynomial algebra $\overline{\mathcal U}_q({\mathfrak n}_\omega)$ with 
generators $z_1,\dots,z_{\ell(w)}$, and  relations 
$$z_{i}z_{j}=q^{(\gamma_i,\gamma_j)}z_{j}z_{i}.$$
They proved that eg. the P.I. degree of ${\mathcal U}_q({\mathfrak n}_\omega)$ 
could be determined from this much simpler algebra in the case where $q=\varepsilon$ is a primitive root of unity. We return briefly to this degree in the last section of this article.

\medskip

The success of the present endeavor rests on the choice of a good basis of 
the associated quasi Laurent algebra  $\overline{\mathcal L}_q({\mathfrak 
n}_\omega)$. While looking at specific cases (\cite{jp})) such a basis was 
found essentially as the ``diagonals'' of the quantized minors of 
Berenstein-Zelevinsky. 

\medskip

Consider the symplectic form $\overline{\mathbb L}$ defined by the relations 
above, that is, the skew symmetric form defined by
$$i<j:\ \overline{\mathbb L}_{i,j}=(\gamma_i,\gamma_j).$$

Recall that a symplectic form may be brought to a block diagonal form by using 
matrices with integer coefficients and determinant $1$.

The center of $\overline{\mathcal L}_q({\mathfrak n}_\omega)$ is given by the 
null space of $\overline{\mathbb L}$.

\medskip

The second major step forward comes with the construction,  quite explicitly, 
and while again using ideas from (\cite{jz},\cite{jp}), of a partial inverse 
$\overline{\mathbb B}$ to $\overline{\mathbb L}$.

\medskip

As an aside, we mention that we actually construct a quantum seed.

\medskip

With these steps taken, the center of $\overline{\mathcal L}_q({\mathfrak 
n}_\omega)$ is easily determined.

\medskip

The final step towards determining the center of ${\mathcal L}_q({\mathfrak 
n}_\omega)$ comes when one realizes that $\overline{\mathbb L}$ actually is 
also the symplectic form ${\mathbb L}$ for a certain family of $q$ commuting 
quantized minors and likewise $\overline{\mathbb B}$ is expressible in terms of 
these minors and their inverses and hence we get a compatible pair  $({\mathbb 
L},{\mathbb B})$.

\medskip

As mentioned in the abstract, the center is given by the null space of 
$$1+\omega.$$

More precisely, to each fundamental weight $\Lambda_s$ there is a covariant 
element $C_s(\omega)$ and the center is given by those $\prod_sC_s^{n_s}(\omega)$ for 
which
$$(1+\omega)(\sum_s n_s\Lambda_s)=0.$$

We discuss some special examples of this in Subsection~\ref{subsec62}, 
Section~\ref{sec7}, and Section~\ref{sec9}.

\medskip

The results have previously been obtained by other methods by Ph. Caldero \cite{Ca0}, \cite{Ca} in the case where $\omega$ is the longest element.

\medskip

Recently, a full treatment was obtained by M. Yakimov \cite{yak} also by different methods.

\medskip

Besides the use of a number of deep results, our methods are completely elementary and focuses on quantum seeds. We also point to certain useful structures related to subalgebras of a fixed parabolic subalgebra ${\mathfrak p}$ and this, we hope, will be seen as significant points of the article. One such structure is the assignment of a unique pair $(c,d)\in{\mathbb N}^2$ to each Schubert Cell. Another is, in the terminology of cluster algebras, a compatible pair given quite explicitly.   

With this available, we can even give a generalization to the centers of algebras connected to double Schubert Cells defined by minimal left coset representatives $\omega^{\mathfrak a}, \omega^{\mathfrak c}$ in $W_{\mathfrak p}\backslash W$ with, say, $\omega^{\mathfrak a}<\omega^{\mathfrak c}$. Here the center is given by the null space of $\omega^{\mathfrak a}+\omega^{\mathfrak c}$

\medskip

Here is a more detailed account of the content:

\smallskip

Section~\ref{sec2}: Background; saturated sets of positive roots are discussed. 
Section~\ref{sec4}: 
It is determined that ${\mathcal U}_q({\mathfrak n}_\omega)$ has the structure of quadratic algebra.
Section~\ref{sec5}: Basics; a diagrammatic way of representing ${\mathfrak 
n}_\omega$ and $\Delta^+({\mathfrak n}_\omega)$ is introduced. 
Section~\ref{sec6}: The case of the associated quasi-polynomial algebra is described and an example is given. 
Section~\ref{sec7} is a Diophantine interlude in which the centers are computed 
for some specific elements $w\in W$ in type $A_n$. Then in Section~\ref{sec8} 
the quantum minors of Berenstein-Zelevinsky are introduced, 
the twist by (\cite{leclerc}) is given, and two series of what we call 
Levendorskii-Soibelman  quadratic algebras are introduced. 
After that, the  way has been  paved for
Section~\ref{sec9} in which the previous results are extended to the general 
setting for these series of quadratic algebras. A more general class of Double Schubert  Cell algebras is also briefly discussed. Finally in the last section, the case where $q$ is a primitive root of unity is discussed.

\section{On Parabolics}

The 
origin of the following lies in A.
Borel \cite{borel}, and B. Kostant \cite{kos}. Other main
contributors are
\cite{bgg} and
\cite{stein}. See also \cite{cap}. We have also found (\cite{sager}) useful.

\medskip

We consider a simple Lie algebra ${\mathfrak g}$.  ${\bf A}$ is the  Cartan 
matrix of ${\mathfrak g}$ and is assumed to be of finite type. $\Pi$ denotes  
a fixed choice of simple roots, and ${\bf E}_{\Pi}$ denotes the euclidean 
space spanned by the simple roots. The fundamental weights corresponding to the 
simple roots are
denoted by
$\Lambda_i$. 

\medskip

\begin{Def} Let $w\in W$. Set $$\Phi_\omega=\{\alpha\in \Delta^+\mid
w^{-1}\alpha\in
\Delta^-\}=w( \Delta^-)\cap  \Delta^+.$$\end{Def}

\medskip

We have that $\ell(w)=\ell(w^{-1})=\vert\Phi_\omega\vert$.

\medskip

\begin{Def}
A subset $S$ of $\Delta^+$ is saturated if whenever 
$\alpha,\beta\in 
S$ and $\alpha+\beta$ is a root, then $\alpha+\beta\in S$.
\end{Def}
\label{sec2}
\begin{Thm}[\cite{kos}]
The map $$w\mapsto \Phi_\omega$$ defines a bijection between $W$ and
the set of all
subsets  $\Phi\subseteq  \Delta^+$ for which both $\Phi$ and $ 
\Delta^+\setminus\Phi$ are saturated. 
\end{Thm}
\medskip

In passing we observe that, trivially, for a saturated set
$\Phi$, both $\Phi$
and $ \Delta^+\setminus\Phi$ correspond to nilpotent
subalgebras.

We will from now on set
$\Phi_\omega=\Delta^+(w)$. 

\medskip

We will consider nilpotent quantized enveloping algebras of the form 
$${\mathcal U}_q(\frak n_\omega),$$
where $\omega$ is an arbitrary element in the Weyl group $W$, and $\frak 
n_\omega$ is the nilpotent defined by the roots 
$\alpha\in\Phi_\omega$ (This is  $\frak 
n_{\omega^{-1}}$ of (\cite{leclerc}). It is convenient for us, also with an eye to 
forthcoming investigations, to assume that we are working with a fixed 
parabolic 
$\mathfrak p$ with a Levi decomposition
\begin{equation}\label{1}
{\mathfrak p}={\mathfrak l}+{\mathfrak u},
\end{equation}
where ${\mathfrak l}$ is the Levi subalgebra, and such that, on the classical 
level, $\frak 
n_\omega\subseteq \frak u$. There is no loss of generality in that.

\medskip

Finally set
\begin{Def}
\begin{eqnarray*}
W_{\mathfrak p}&=&\{w\in W\mid \Phi_\omega\subseteq \Delta^+({\mathfrak l})\}\\
W^{\mathfrak p}&=&\{w\in W\mid \Phi_\omega\subseteq \Delta^+({\mathfrak u})\}.
\end{eqnarray*}
$W^{\mathfrak p}$ is a set of distinguished representatives of the right
coset space
$W_{\mathfrak p}\backslash W$.
\end{Def}

\medskip

It is well known (see eg (\cite{sager})) that any $w\in W$ can be written 
uniquely as $w=w_pw^{p}$ with $w_p\in W_{\mathfrak p}$ and $w^{p}\in W^{\mathfrak p}$.

\section{The quadratic algebras}\label{sec4}

Let $\omega=s_{\alpha_1}s_{\alpha_2}\dots
s_{\alpha_t}$ be an element of the Weyl group written in reduced form. 
Using his  braid operators, given in a special
case as
$$T_i(e_j)=\sum_{a+b=r}(-q_i)^{-b}e_i^{(a)}e_je_i^{(b)},$$ where $r=-\langle  
h_i,\alpha_j\rangle$, Lusztig in (\cite{luz}) construct a sequence of elements 
$Z_1,\dots,Z_t \subseteq {\mathcal U}^+_q({\mathfrak g})$. Specifically, 
$$Z_{i}=T_{\omega_{i-1}}(E_{\alpha_{i}}),\ i=2,\dots,t,\textrm{
and
}Z_1=E_{\alpha_1},$$
where, for each $i=1,\dots,t$,  
$\omega_i=s_{\alpha_1}s_{\alpha_2}\dots s_{\alpha_{i}}$. In
particular,
$\omega=\omega_t$. The weight of $Z_{i}$ is given by 
$\gamma_i=\omega_{i-1}(\alpha_i)$. Here, and throughout, we use the notation of (\cite{jan}).

The following result is well known

\begin{Thm}[\cite{lev},\cite{lev0}] Suppose that  $1\leq i<j\leq t$. Then
$$Z_{i}Z_{j}=q^{( 
\gamma_i,\gamma_j)}Z_{j}Z_{i} +
\textrm{terms involving only elements $Z_k$ with $i\leq k\leq
j$.}$$
\end{Thm}
Our statement follows \cite{jan},\cite{jan2}. Other authors, eg.
\cite{lev},
\cite{leclerc} have used the other Lusztig braid operators. The
result is
just a difference between $q$ and $q^{-1}$. Proofs of this
theorem which are
more accessible are available (\cite{cp},\cite{jan2}).

\medskip

\begin{Prop}
${\mathcal U}_q({\mathfrak n}_\omega)$ is a quadratic algebra.
\end{Prop}

It is known that this algebra is isomorphic to the algebra of functions on 
${\mathcal U}_q({\mathfrak n}_\omega)$ satisfying the usual finiteness 
condition. It is analogously equivalent to the algebra of functions on 
${\mathcal U}^-_q({\mathfrak n}_\omega)$ satisfying a  similar finiteness 
condition. See eg (\cite{leclerc}) and (\cite{jan}). We will not distinguish 
between these algebras.
\bigskip

\section{basic structure}\label{sec5}

Consider a fixed basis 
$\Pi=\{\alpha_1,\alpha_2,\dots,\alpha_R\}$ of $\Phi$. Let us agree to write 
$\sigma_{\alpha_i}$ in $W$ just as $\sigma_i$ for $i=1,\dots,R$. Returning to the decomposition (\ref{1}), assume henceforth that  ${\mathfrak u}\neq\{0\}$.

\smallskip

Adapting to the language of \cite{fz}, \cite{bz}, and 
others, we will often label structures derived from  $\omega^{\mathfrak 
r}$ by the reduced word  ${\mathfrak 
r}$. The full structure with double words will not be required here.

\medskip

\medskip

Let $\omega^{\mathfrak p}$ be the maximal element in $W^{\mathfrak p}$. It
is the one which
maps all roots in $\Delta^+({\mathfrak
u})$ to $\Delta^-$. (Indeed: To $\Delta^-({\mathfrak u})$.)
Let $\omega_0$ be the
longest element in $W$ and $\omega_L$ the
longest in the Weyl group of ${\mathfrak l}$, Then
\begin{equation}\omega^{\mathfrak p}\omega_L=\omega_0.\end{equation}

Let 
$\omega^{\mathfrak 
r}=\sigma_{i_1}\sigma_{i_2}\cdots\sigma_{i_r}\in W^{\mathfrak p}$ be fixed  and written in 
a fixed reduced form. Then
$\ell(\omega^{\mathfrak 
r})=r$. 

 Set  
\begin{equation}\Delta^+(\omega^{\mathfrak 
r})=\{ \beta_{i_1},\dots,\beta_{i_r}\}.\end{equation}

Consider $\alpha_{i_{r+1}}\in\Pi$ and set $\omega_1=\omega^{\mathfrak 
r}\circ
s_{i_{r+1}}$.

\noindent{\bf Case 1:}
$\omega^{\mathfrak 
r}(\alpha_{i_{r+1}})=\gamma\in\Delta^+(\mathfrak u)$.
Then $\ell(w_1)=\ell(\omega^{\mathfrak 
r})+1$ and
\begin{equation}\Delta^+(w_1)= \{
\beta_{i_1},\dots,\beta_{i_r}\}\cup\{\gamma\}\textrm{ and
}\gamma=\beta_{i_{r+1}}=\omega^{\mathfrak 
r}(\alpha_{i_{r+1}}).\end{equation}
Thus, $\omega_1\in W^{\mathfrak p}$.

\smallskip

\noindent{\bf Case 2:}
$\omega^{\mathfrak 
r}(\alpha_{i_{r+1}})=-\gamma\in\Delta^-(\mathfrak u)$.
Then $\ell(w_1)=r-1$ and
\begin{equation}\Delta^+(w_1)= \{
\beta_{i_1},\dots,\beta_{i_r}\}\setminus\{\gamma\}\textrm{
and }\gamma=\beta_{i_s}\textrm{ for some
}\beta_{i_s}\in\Delta^+(\omega^{\mathfrak 
r}).\end{equation}

\medskip

We must always be in at least one of these cases since otherwise $\omega^{\mathfrak r}$ would map all simple roots, hence $\Delta({\mathfrak g})$, to $\Delta({\mathfrak l})$. If $\omega^{\mathfrak 
r}$ is maximal and $\alpha_i$ is a simple root such that  $\omega^{\mathfrak 
r}(\alpha_i)\in\Delta^-({\mathfrak u})$ then $\alpha_i\in \Delta^+({\mathfrak u})$. It is easy to see that under the same assumptions, $\omega^{\mathfrak r}(\alpha_i)\in\Delta^-({\mathfrak l})$ is not possible. 
Furthermore, if $\alpha_i\in\Delta^+({\mathfrak l})$ then  $\omega^{\mathfrak r}(\alpha_i)\in\Delta^+({\mathfrak g})$. In conclusion, a maximal $\omega^{\mathfrak r}$ maps $\Delta^+({\mathfrak l})$ to $\Delta^+({\mathfrak l})$
and $\Delta^+({\mathfrak u})$ to $\Delta^-({\mathfrak u})$. Thus, if $\omega^{\mathfrak r}$ is maximal, $\omega_L\omega^{\mathfrak r}=\omega^{\mathfrak r}\omega_L=\omega_0$. Hence $\omega^{\mathfrak r}=\omega^{\mathfrak p}$.

\medskip

It follows easily that we have the following
conclusion: Let $\omega^{\mathfrak 
r}\in W^{\mathfrak p}$ with $\ell(\omega^{\mathfrak 
r})=r$. Then we may write \begin{equation}\omega^{\mathfrak 
r}=s_{i_1}\circ \dots \circ s_{i_s}\circ\dots\circ
s_{i_r}\textrm{ where for all
}j=1\dots,r: w_j=s_{i_1}\circ \dots
\circ s_{i_j}\in W^{\mathfrak p}.\end{equation}Furthermore,
\begin{equation}\Delta^+(\omega^{\mathfrak 
r})= \{
\beta_{i_1},\dots,\beta_{i_s},\dots,\beta_{i_r}\}\end{equation}where for 
each $i_s$:\begin{equation}\beta_{i_s}=s_{i_1}\circ
\dots \circ s_{i_{s-1}}(\alpha_{i_s})\textrm{ for $s>1$, and 
}\beta_{i_1}=\alpha_{i_1}.\end{equation}
Moreover,

\begin{equation}\label{eq10}\omega^{\mathfrak 
r}=\sigma_{{i_1}}\sigma_{{i_2}}\cdots 
\sigma_{{i_r}}=\sigma_{\beta_{r}}\cdots 
\sigma_{\beta_{2}}\sigma_{\beta_{1}}.\end{equation}

\medskip

From now on,  $q$ is a fixed element of
${\mathbb C}$ which
is 
not 
a root of unity, $\omega^{\mathfrak r}\in W^{\mathfrak p}$ is given with a fixed decomposition as in (\ref{eq10}), and $\Delta^+(w^{\mathfrak r})$ is
our universe.

\begin{Def}
Let
${\mathbf b}$
denote the
map $\Pi\to\{1,2,\dots,R\}$ defined by ${\mathbf
b}(\alpha_i)=i$.
Let  $\overline\pi_{\mathfrak r}:\{1,2,\dots,
r\}\to\Pi$ be given by
\begin{equation}\overline\pi_{\mathfrak
r}(j)=\alpha_{i_j}.\end{equation}
If $\overline\pi_{\mathfrak r}(j)=\alpha$ we say that $\alpha$
(or
$\sigma_\alpha$) occurs at position 
$j$ 
in $\omega^{\mathfrak r}$, and we say that
$\overline\pi_{\mathfrak r}^{-1}(\alpha)$ are the positions at
which $\alpha$
occurs in $w$.
Set \begin{equation}
{\pi}_{\mathfrak r}={\mathbf b}\circ\overline\pi_{\mathfrak r}.       
\end{equation}
\end{Def}

\bigskip

Let, for $1\leq n\leq r, \ 
\omega_n=\sigma_{\alpha_{i_1}}\sigma_{\alpha_{i_2}}\cdots 
\sigma_{\alpha_{i_{n}}}$. Thus, we have a 1-dimensional 
presentation of the situation given by the (ordered) set
$\{1,2,\dots,r\}$.

The following 2-dimensional presentation is even more useful and
informative:

\medskip

\begin{Def}
\begin{eqnarray}&{\mathbb U}({\mathfrak
r})=\\&\{(s,t)\in{\mathbb
N}\times {\mathbb N}\mid 
\exists n \textrm{ 
such that } s=\pi_{\mathfrak r}(n)\textrm{ and }
\omega_n=\omega_1
\sigma_{i_n}\omega_2\dots\omega_t
\sigma_{i_n}\}.\nonumber\end{eqnarray}
\end{Def}
In the above, it is understood that each $\omega_i\in
W\setminus\{e\}$ is
reduced and does 
not contain any $\sigma_{i_n}$.

We also identify $n\leftrightarrow(s,t)$ (and  
$\beta_{n}\leftrightarrow\beta_{s,t}$) if $n,s,t$ are
connected as
above. We denote the identifications between $\{1,2,\dots,r\}$
and ${\mathbb
U}({\mathfrak r})$ simply as $n\leftrightarrow (s,t)$. (And 
$\beta_{n}\leftrightarrow \beta_{s,t}$)

We define a map $\pi_{\omega_n}$ for such $\omega_n$ in analogy
with that of
$\pi_{\mathfrak r}$.

If $\omega^{\mathfrak r}=\omega_m\tilde\omega$ and 
$\omega_{m}=\omega_n\hat\omega$ with $\omega_n,\omega_m\in W^{\mathfrak p}$
and all Weyl
group elements reduced, we say that $\omega_n<\omega_m$ if $\hat\omega\neq e$. 
\begin{Def}If $n\leftrightarrow 
(s,t)$ and $m\leftrightarrow (c,d)$ we define
\begin{equation}(s,t)<(c,d)
\Leftrightarrow \omega_{s,t}<\omega_{c,d}.\end{equation}
\end{Def}

\medskip

For a fixed $s\in\{1,2,\dots,R\}$ we let $s_{\mathfrak r}$ denote the
maximal such $t$.
This is the number of tines $\sigma_s$ occurs in $\omega^{\mathfrak
r}$. We then have
\begin{equation}
{\mathbb U}({\mathfrak r})=\{(s,t)\in {\mathbb N}\times {\mathbb
N}\mid
1\leq s\leq 
R\textrm{ 
and }1\leq t\leq s_{\mathfrak r}\}.
\end{equation}

\medskip

Finally, notice that if $(s,t)\in{\mathbb U}({\mathfrak r})$ then we may
construct a subset
${\mathbb U}({\mathbf s}, {\mathbf t})$ of ${\mathbb U}$ by the above recipe,
replacing
$\omega^{\mathfrak r}$ by $\omega_{s,t}$. In this subset $t$ is
maximal.

\bigskip

\section{The quasi-polynomial algebra}
\label{sec6}
\subsection{The first definitions and computations}
\begin{Def} \label{4.22}
$\overline{\mathcal M}_q({\mathfrak r})$ denotes the ${\mathbb C}$-algebra
generated by elements 
$\{z_{j};j=1,\dots,r\}$ indexed by the elements $\beta_{j}$
defined as above
and 
with relations
\begin{equation}z_iz_j=q^{(\beta_{i},\beta_{j})}z_jz_i
\textrm{ if
} 
i<j.\end{equation}
We let $\overline{\mathcal L}_q({\mathfrak r})$
denote the
associated quasi-Laurent algebra, and we let $\overline{\mathcal
X}_q({\mathfrak
r})$ denote the center of $\overline{\mathcal L}_q({\mathfrak
r})$. We will also label the generators by $z_{s,t}$ as discussed in 
Section~\ref{sec5}.
\end{Def}
\bigskip

Let ${s}\in Im(\pi_{\mathfrak r})$. It is then straightforward
to see that
\begin{equation}
 -\omega(\Lambda_s)+\Lambda_{s}=\beta_{s,1}+\beta_{s,2}
+\dots+\beta_{s,s_{\mathfrak r}}.
\end{equation}

\begin{Def}Let ${s}\in Im(\pi_{\mathfrak r})$.
We define the  element $\overline 
C_{s}({\mathfrak r})$ in the 
quasi-polynomial algebra $\qpol$ by 
\begin{equation}\label{eq1}\overline C_{s}({\mathfrak 
r})=z_{s,1}z_{s,2}\cdot\dots\cdot 
z_{s,s_{\mathfrak r}}.\end{equation}
\end{Def}

\medskip

\begin{Prop}\label{commu1}The following holds for all $s\in
Im(\pi_{\mathfrak
r})$ and all 
$(a,b)\in{\mathbb U}({\mathfrak r})$:
\begin{equation}z_{a,b} \overline C_{s}({\mathfrak
r})=q^{-(
\beta_{a,b},\left(1+\omega^{\mathfrak
r}\right)(\Lambda_{s}))}
\overline C_{s}({\mathfrak r})z_{a,b},\end{equation}
\end{Prop}

\proof We will be using repeatedly that
$\beta_{s,t}=-\omega_{s,t}(\alpha_{i_s})$. Consider a decomposition $\omega^{\mathfrak 
r}=\omega_A\sigma_{\alpha_{i_\ell}}\omega_B$ with
$\pi_{\mathfrak
r}(\ell)=a\neq s$ and suppose that
$\beta_{s,t}<\beta_{a,b}<\beta_{s,t+1}$ for
some $t$. It is here understood that $\ell=(a,b)$. The elements 
$\omega_A,\omega_B$ of course depend on $\ell$, indeed, 
$\omega_A\sigma_{i_\ell}=\omega_\ell$.
Then
\begin{eqnarray}\nonumber z_\ell z_{s,t+1}\dots
z_{s,s_{\mathfrak r}}&=&q^{(
\beta_{i_\ell},(-\omega^{\mathfrak r}+\omega_{A})
(\Lambda_{s})
))} 
z_{s,t+1}z_{s,t+2}\dots z_{s,s_{\mathfrak r}}z_\ell\\&=&q^{(
\alpha_{i_\ell},(-\sigma_{\alpha_{i_\ell}}\omega_{B}
(\Lambda_{s})+\Lambda_{s}
))} 
z_{s,t+1}z_{s,t+2}\dots z_{s,s_{\mathfrak r}}z_\ell\\
&=&q^{(
\alpha_{i_\ell},\omega_{B}(\Lambda_{s}))}z_{s,t+1}\dots
z_{s,s_{\mathfrak r}}z_\ell.\end{eqnarray}

Similarly, 

\begin{equation}
 z_\ell  z_{s,1}\dots z_{s,t}=q^{-( 
{\alpha_{i_\ell}},(\omega^{-1}_A(\Lambda_{s}))}
z_{s,1}\dots
z_{s,{t}}z_\ell .
\end{equation}
  
The statement then follows directly.   To complete this part of the picture, we need to consider $z_\ell<z_{s,1}$
and
$z_\ell>z_{s,s_{\mathfrak r}}$

The case $z_\ell<z_{s,1}$ easily results in the exponent
${(
\beta_\ell,\left(1-\omega^{\mathfrak
r}\right)(\Lambda_{s}))}$, but here
${( 
\beta_\ell,\Lambda_{s})}=0$. The case
$z_\ell>z_{s,s_{\mathfrak r}}$ gives an
exponent ${-( 
\beta_\ell,\left(1-\omega^{\mathfrak
r}\right)(\Lambda_{s}))}$, and here
$( \beta_\ell,\omega^{\mathfrak r}(\Lambda_{s})=0$.

Next, we observe the following simple formulas, where $\omega_A$ and $\omega_B$ now are determined by $(s,t)$:
\begin{equation} z_{s,t} z_{s,t+1}\dots
z_{s,s_{\mathfrak r}}=q^{-1+(
\alpha_{s},\omega_{B}(\Lambda_{s}))}  z_{s,t+1}\dots 
z_{s,s_{\mathfrak r}}z_{s,t},\end{equation}
and, similarly,
\begin{equation} z_{s,t}z_{s,1}\dots z_{s,t-1}=q^{1-( 
\alpha_{s},\omega_A^{-1}(\Lambda_{s}))} z_{s,1}\dots 
z_{s,t-1}z_{s,t}.\end{equation}

These formulas also hold at the extreme positions of $z_{s,1}$
and
$z_{s,s_{\mathfrak r}}$, where either $\omega_A=1$ or $\omega_B=1$.

So,   indeed for any simple root $\alpha=\pi_\omega(\ell)$ and decomposition 
$\omega^{\mathfrak
r}=\omega_\ell\omega_B=\omega_As_\alpha\omega_B$, we get,
for the corresponding $z_\ell$,

\begin{equation}
 z_\ell \overline C_{s}({\mathfrak r})=q^{( 
\alpha,\left(\omega_B-\omega_A^{-1}\right)(\Lambda_{s}))}\overline 
C_{s}({\mathfrak r})z_\ell,
\end{equation}
which, by the previous definitions is equivalent to the statement in the 
proposition.\qed

If $s\notin Im(\pi_{\mathfrak r})$, we set $\overline
C_{s}({\mathfrak r})=1$.
To any linear combination
$$\sum_in_i\Lambda_{i}$$ with integer
coefficients we may consider the element 
in the quasi-Laurent algebra
\begin{equation}
 (\overline C_{s_1}({\mathfrak r}))^{n_1}\dots (\overline 
C_{s_k}({\mathfrak r}))^{n_k}.
\end{equation}
 If $\overline C_{s}({\mathfrak 
r})=1$  we set $n_s=0$.

\medskip

Let 
$$S_{\mathfrak r}=\textrm{Span}\{\alpha_s\mid s\in
Im(\pi_{\mathfrak r})\}.$$
It is obviously invariant under $\omega^{\mathfrak r}$. 
We view tacitly the elements $\Lambda_s$ as restricted to this
space.

\medskip

\begin{Prop}\label{usef}Let $n_1,\dots,n_k$ be integers and let 
$s_1,\dots,s_k\in Im(\pi_{\mathfrak r})$.
$$(\overline C_{s_1}({\mathfrak r}))^{n_1}\dots (\overline 
C_{s_k}({\mathfrak r}))^{n_k}\in \overline{\mathcal
Z}_q({\mathfrak
r})\Leftrightarrow
(1+\omega^{\mathfrak r})(\sum_{j=1}^k n_j\Lambda_{s_j})=0 .$$
\end{Prop}
 
 \proof
 The 
commutation between this and any 
$z_\ell$ is given by
\begin{eqnarray}\nonumber  z_\ell  (\overline C_{s_1}({\mathfrak 
r}))^{n_1}\dots (\overline 
C_{s_k}({\mathfrak r}))^{n_k}= \\  q^{-( 
\beta_\ell,\left(1+\omega^{\mathfrak r})\right)(\sum_i 
n_i\Lambda_i))} (\overline C_{s_1}({\mathfrak
r}))^{n_1}\dots
(\overline 
C_{s_k}({\mathfrak r}))^{n_k} z_\ell. \end{eqnarray}\qed
 
 \medskip
 
This actually determines the center as will be proved below after some
preparation.

\begin{Rem}
Since $1+\omega^{\mathfrak r} $ is an integer matrix, there is an ${\mathbb R}$ 
basis of the null space given by vectors with integer coordinates in the basis 
of fundamental weights. 
\end{Rem}

\bigskip

\subsection{More definitions and computations. The center of the \label{subsec6}
quasi-polynomial algebra. }

\flushleft 

We first make a very useful observation:

\begin{Lem}\label{3.1}  Let $\alpha_i\in \Phi$. Then 
$$(s_i+1)(\Lambda_i)+\sum_{j\neq 
i}a_{ji}(\Lambda_j)=0.$$
\end{Lem}

\proof This is actually equivalent to \cite[(2.27)]{fz}
by the definition of  $a_{ki}$. \qed

\bigskip

\begin{Def}Let $(s,t)\in{\mathbb U}({\mathfrak r})$.
Set\begin{equation} 
 {\overline M}^\downarrow_{s,t}=z_{s,1}\dots z_{s,t}.\end{equation}
\end{Def}

This element has weight 
\begin{equation}
 {p_{s,t}=-w_{s,t}(\Lambda_{s})+\Lambda_{s}}=\beta_{s,1}+\dots+ 
\beta_{s,t}.\label{28}
\end{equation}

\begin{Prop}\label{5.9}Let $(a,b),(s,t)\in{\mathbb U}({\mathfrak r})$. Then 
\begin{equation}
z_{a,b}\overline M^{\downarrow}_{s,t}=q^{E_{a,b}} \overline 
M^{\downarrow}_{s,t}z_{a,b}, 
\end{equation}where the exponent $E_{a,b}$
is given as follows:
$CASE\ 1: (a,b)\leq (s,t) $:  
\begin{equation}
 E=-(\beta_{a,b},(1+w_{s,t})(\Lambda_{s})).   
\end{equation}
$CASE\ 2: (a,b)> (s,t)$:
\begin{equation}
 E=-(\beta_{a,b},(1-w_{s,t})(\Lambda_{s})).  
\end{equation}
\end{Prop}
\proof CASE 1 is equivalent to Proposition~\ref{commu1} and CASE
2 follows by
very similar arguments.\qed

As a special case we get (because we here are in Case 1 only)
\begin{Cor}\label{only}
\begin{equation}\forall (s,t)\in {\mathbb U}({\mathfrak r}),\forall j\in 
Im(\pi_{\mathfrak r}): \overline
M^{\downarrow}_{(s,t)}\overline
C_{j}({\mathfrak r})=q^{-( 
(1-\omega_{s,t}(\Lambda_{s}),(1+\omega^{\mathfrak 
r})(\Lambda_{j}))}\overline 
C_{j}\overline M^{\downarrow}_{s,t}.\end{equation}
\end{Cor}

\medskip

The following formula, which is only seemingly more general, and  in which $(a,b),(c,d)\in {\mathbb 
U}({\mathfrak 
r})$,  is also 
useful. 
\begin{Cor}\label{5.11}\begin{equation} \label{5.11-f} 1\leq 
(s,t){\leq}(c,d)\Rightarrow\overline M_{s,t}^\downarrow \overline 
M_{c,d}^\downarrow=q^{-((1-\omega_{s,t})(\Lambda_{s}),(1+\omega_{c,d})(\Lambda_{
c}))}
\overline M_{c,d}^\downarrow\overline
M_{s,t}^\downarrow.\end{equation}
\end{Cor}

\bigskip

\medskip

\begin{Def}
We set $\omega_{s,0}(\Lambda_s)=\Lambda_s$ and $\overline M_{s,0}^\downarrow=1$ 
for all $s\in Im(\pi_{\mathfrak r})$.
\end{Def}

Notice that $\omega_{s,s_{\mathfrak r}}(\Lambda_s)=\omega^{\mathfrak
r}(\Lambda_{s})$.

\bigskip

For a fixed, but arbitrary, $(s,t)\in{\mathbb U}({\mathfrak r})$
and $a\in
Im(\pi_{\omega_{s,t}})$ set
$\overline{p}_{\mathfrak r}(a,s,t)=\max(\pi_{\omega_{s,t}}^{-1}(a))$. If 
$a\notin
Im(\pi_{\omega_{s,t}})$ we set $\overline{p}_{\mathfrak r}(a,s,t)=0$.

\medskip

Then we define
\begin{Def}\label{6.12}
\begin{eqnarray}\overline{F}(s,t)&=&
\overline M^\downarrow_{s,t}\overline
M^\downarrow_{s,t-1}\prod_{a_{js}<0}
\left(\overline
M^\downarrow_{j,\overline{p}_{\mathfrak 
r}(j,s,t)}\right)^{a_{js}}.\end{eqnarray}
\end{Def}
Recall that ${\mathfrak g}$ is finite-dimensional, hence diagonalizable. It is then a standard fact that $d_s=d_{\alpha_s}=(\Lambda_s,\alpha_s)$.  The key result is: 

\begin{Prop}\label{key}The following holds holds for all $(c,d),(s,t)\in{\mathbb
U}({\mathfrak
r})$:
\begin{equation}
 z_{{c,d}}\overline{F}(s,t)=q^{E_{c,d}} \overline{F}(s,t)z_{{c,d}},
\end{equation}where
\begin{equation}
 E_{c,d}=-( \beta_{c,d},\alpha_{s}) 
+2{d_s}\delta_{(c,d),(s,t)}.
\end{equation}
Furthermore,
\begin{equation}
 \overline M^\downarrow_{{c,d}}\overline{F}(s,t)=q^{G_{c,d}} 
\overline{F}(s,t)\overline M^\downarrow_{{c,d}},
\end{equation}where
\begin{equation}
 G_{c,d}=-( (\beta_{c,1}+\dots +\beta_{c,d}),\alpha_{s}) 
+2{d_s}(\delta_{(c,1),(s,t)}+\dots+\delta_{(c,d),(s,t)}) . 
\end{equation}
\end{Prop}

\proof  Our strategy is to sum up the exponents of the the $q$-commutation for 
the minors of weights 
$-\omega_{s,t}\circ\sigma_{s}(\Lambda_{s})+\Lambda_{s}$, 
$-\omega_{s,t}(\Lambda_{s})+\Lambda_{s}$ and 
$a_{js}(-\omega_{s,t}(\Lambda_{j})+\Lambda_{j}), j\neq s$  and then use 
the formulas of Proposition~\ref{5.9} in combination with
Lemma~\ref{3.1}. This goes well and yields exponent
$E=-(\beta_{c,d},\alpha_{s})$ in all cases where
$(c,d)\neq(s,t)$ since there
is no mixing of Case 1 and Case 2 positions for the various
minors. When $(c,d)=(s,t)$  we do mix Case 1 and Case 2, but only in
reference to the minor of weight
$-\omega_{s,t}\circ\sigma_{s}(\Lambda_{s})+\Lambda_{s}$.
Here we notice  that $2\Lambda_{s}+\sum_{j\neq
s}a_{js}\Lambda_j=
(1-\sigma_{s})\Lambda_{s}=\alpha_{s}$ by the same 
lemma. The results follow. \qed

\bigskip

Now we combine two (adjacent) minors based on the same
$\alpha_{\pi_\omega(n)}$ but of
opposite signs. In the following we have to assume that $t<s_{\mathfrak r}$ since we 
work with both $t$ and $t+1$.

The following follows then easily:

\begin{Lem}
\begin{eqnarray}\nonumber &\forall (c,d)\in {\mathbb U}({\mathfrak r}), \forall 
(s,t)\in {\mathbb U}({\mathfrak r}) \textrm{ satisfying }t<s_{\mathfrak r},:\\&
z_{c,d}
\overline{F}(s,t)\overline{F}(s,t+1)^{-1}=q^{\nabla_{c,d}}\overline{F}(s,
t)\overline{F}(s,t+1)^{-1}z_{c,d}\end{eqnarray}
where
\begin{equation}\nabla_{c,d}=2{d_s}\delta_{(c,d),(s,t)}-2{d_s}\delta_{(c,d),(s,t+1)}
.\end{equation}
\end{Lem}

We can thus construct some elements that commute with 
everything except a single element:
\begin{Def}$\forall 
(s,t)\in {\mathbb U}({\mathfrak r})$ satisfying $t<s_{\mathfrak r}$
set\begin{equation}\label{6.15}\overline{B}_{s,t}=\overline{F}(s,t)\overline{F}
(s,t+1)^{-1}
.\end{equation}\end{Def}

\begin{Prop}\label{6.16}Whenever $\overline{B}_{s,t}$ is defined, the following 
holds for 
all $(c,d)\in{\mathbb
U}({\mathfrak
r})$:
 $$\overline M^\downarrow_{c,d} 
\overline{B}_{s,t}=q^{2{d_s}\delta_{(c,d),(s,t)}}\overline{B}_{s,t}\overline M
^\downarrow_{c,d}.$$
\end{Prop}

\begin{Def}We define a symplectic form $\overline{\mathbb L}_0$ by, for 
$\beta_i,\beta_j\in \Delta^+(\omega^{\mathfrak r})$
\begin{equation}
 \beta_i<\beta_j\Rightarrow \overline{\mathbb L}_{0,ij}=(\beta_{i},\beta_{j})
\end{equation}
We let $\overline{\mathbb L}$ denote the symplectic form defined in terms of the 
elements  $\overline M^\downarrow_{c,d}$;
\begin{equation}
\overline M^\downarrow_{c,d}\overline M^\downarrow_{s,t}=q^{{\mathbb 
L}_{(c,d),(s,t)}}\overline M^\downarrow_{s,t}\overline M^\downarrow_{c,d}.
\end{equation}
We will also find it convenient to introduce an auxiliary sesquilinear form in 
which the elements $\overline M^\downarrow_{c,d}$ form an orthonormal basis.
\end{Def}

We define a matrix ${\mathbb A}$ by
\begin{equation}
{\mathbb A}_{(c,d),(s,t)}=\left\{\begin{array}{l}1\textrm{ if }s=c, 
t=1,2,\dots,d\\0\textrm{ else }\end{array}\right. .
\end{equation}

It follows from (\ref{28}) that

\begin{Lem}
$\overline{\mathbb L}$ is obtained, as a form, by a change of basis by the 
formula
\begin{equation}
\overline{\mathbb L}={\mathbb A}\overline{\mathbb L}_0{\mathbb A}^t.
\end{equation}
\end{Lem}

\medskip

It now follows from Proposition~\ref{6.16} that we have

\medskip
\begin{Cor}\label{cor}For any $(s,t)\in{\mathbb U}({\mathfrak
r})$ for which
$t<s_{\mathfrak r}$, 
the (basis) vector $2{d_s}\overline M^\downarrow_{s,t}$ belongs to
the image of
$\overline{\mathbb L}$. 
\end{Cor}

\medskip

Since

\begin{equation}\textrm{Range}(\overline{\mathbb 
L})^\perp=Ker(\overline{\mathbb 
L})\end{equation}

it follows that the kernel is contained in the space spanned by
the elements
$\overline C_s({\mathfrak r})$. Since $q$ is not a root of unity, the center of 
$\qpol$ is contained in the kernel of $\overline{\mathbb L}$. Combining 
Corollary~\ref{cor}
with
Proposition~\ref{usef}, we have proved the following:
\begin{Thm}[\cite{BCL}, Theorem 3.1; \cite{yak2}, Theorems 4.1 and 4.5]The center of $\overline{\mathcal L}_q({\mathfrak 
r})$ is given by the kernel of 
$(1+\omega^{\mathfrak r})$ on $S_{\mathfrak r}$.
\end{Thm}

\proof We need only consider the eigenspaces of $\omega^{\mathfrak r}$ (in the 
span of the covariant elements), and here it is only the $\pm1$ eigenspaces 
that merit attention: Since we are working with matrices with integer 
coefficients, we see that the mentioned eigenspaces are spanned by elements 
with 
integer coefficients. Observe also that $1+\omega^{\mathfrak r}=2-(
1-\omega^{\mathfrak r})$. 
Consider $\overline C_{+}({\mathfrak r})=\prod_{i\in\textrm{Im}({\pi}_{\mathfrak
r})}(\overline C_{i}({\mathfrak r}))^{n_i} \in \overline{\mathcal
M}_q({\mathfrak
r})$ for which 
\begin{equation}(1-\omega^{\mathfrak r})(\sum_{i\in\textrm{Im}({\pi}_{\mathfrak
r})}
n_i\Lambda_i)=0 .\end{equation}
Then, by
Proposition~{\ref{commu1}, for all $(a,b)$,
\begin{equation} z_{a,b} \overline C_{+}({\mathfrak 
r})=q^{2\beta_{a,b},(\sum_{i\in\textrm{Im}({\pi}_{\mathfrak
r})}
n_i\Lambda_i))} 
\overline C_{+}({\mathfrak r})z_{a,b}.
\end{equation}
Recall that $q$ is not a root of unity. It follows that $\overline 
C_{+}({\mathfrak r})$ can only commute with all $z_{a,b}$ if 
$\sum_{i\in\textrm{Im}({\pi}_{\mathfrak
r})}
n_i\Lambda_i=0$.
\qed

\medskip

\begin{Rem}The kernel of $1+\omega^{\mathfrak r}$ is of course unchanged if we 
enlarge $S_{\mathfrak r}$ to  ${\bf E}_{\Pi}$  by the elements $\Lambda_i$ 
that are left invariant by $\omega^{\mathfrak r}$.\label{ofcourse}
\end{Rem}

\bigskip

\subsection{Example: The full nilpotent - the longest element $\omega_0\in 
W$}\label{subsec62}

 If $\omega=\omega_0$ is the longest word in the Weyl group,
we know that it
is either $-1$ or implemented by a diagram symmetry of order 2.
Indeed, we get $\omega_0=-1$ in all cases except 
\begin{equation}\label{exc1}A_n, D_{2n+1}(n>1), \textrm{ and }E_6.\end{equation}
(See eg. exercises 18, 32 in Chapter 2 in (\cite{knapp}).)
\begin{Prop}
If the simple Lie algebra is not in the list (\ref{exc1}), the center of 
$\qpol=\overline{\mathcal U}_q({\mathfrak n})$ is generated by the elements 
$\overline 
C_{s}({\mathfrak r})$, with $s=1\dots,R$. In the following we use the numbering of simple roots from (\cite{hu}).
In type ${\mathbf A}_\ell$, with simple roots
$\alpha_1,\alpha_2,\dots,
\alpha_\ell$ we have that 
\begin{equation}\omega_m(\Lambda_{\alpha_s})=-\Lambda_{\alpha_{\ell-s+1}},\end{equation}
and the center is generated by the elements\begin{eqnarray*}&\overline 
C_{s}({\mathfrak r})\overline C_{\ell-s+1}({\mathfrak r})\textrm{ for 
}2s\neq\ell+1\\
&\overline C_{s}({\mathfrak r})\textrm{ for }2s=\ell+1\end{eqnarray*}
In type ${\mathbf D}_{2\ell+1}$, the center is generated by the 
elements\begin{equation}\overline C_{s}({\mathfrak 
r})(s=1,\dots,2\ell-1),\textrm{ and }\overline C_{2\ell}({\mathfrak 
r})\overline C_{2\ell+1}({\mathfrak r}).\end{equation}
In type ${\mathbf E}_6$, the center is generated by the elements 
\begin{equation}
\overline C_{1}({\mathfrak r})\overline C_{6}({\mathfrak r}), \overline 
C_{2}({\mathfrak r})\overline C_{5}({\mathfrak r}), \overline C_{2}({\mathfrak 
r}),\textrm{ and }\overline C_{4}({\mathfrak r}).\end{equation}
\end{Prop}

\bigskip

\section{Diophantine interlude}\label{sec7}
Consider a simple Lie algebra of type $A_{a+b+c-1}$ and let
$N=a+b+c$. The
simple roots are
\begin{equation}\Pi=\{e_i-e_{i+1}\mid i=1,2,\dots,
N-1\}.\end{equation}
We choose a Levi subalgebra defined by
\begin{equation}\Sigma=\Pi\setminus\{(e_a,e_{a+1}),(e_{a+b},e_{a+b+1})\}.\end{equation}

We let $I_n$ denote the $n\times n$ identity matrix. We are
interested in
\begin{equation}\omega^{\mathfrak r}=\left(\begin{array}{ccc}
0&0&I_c\\0&I_b&0\\I_a&0&0\end{array}
\right),\end{equation} but introduce a more general family 
of $(a+b+c)\times(a+b+c)$ matrices 
\begin{equation}w([a,\varepsilon_a];[b,\varepsilon_b];[c,\varepsilon_c])=I_{
a+b+c}+\left(\begin{array}{ccc}
0&0&\varepsilon_cI_c\\0&\varepsilon_bI_b&0\\\varepsilon_aI_a&0&0\end{array}
\right).
\end{equation}
Here $\varepsilon_d$ denotes an integer, $d=a,b,c$, such that
$\varepsilon_d^2=1$. We will
always start with $\varepsilon_a=\varepsilon_b=\varepsilon_c=1$
but we will
later encounter more general signs when we perform Gaussian
Elimination moves.
In retrospect, it can be seen that all the matrices encountered here  satisfy 
$\varepsilon_a\varepsilon_b\varepsilon_c=1$. Our task is to try
to determine
${\mathrm
corank}
(w([a,\varepsilon_a];[b,\varepsilon_b];[c,\varepsilon_c]))$.

\bigskip

For arbitrary signs we have, by easy Gaussian moves, the
following reductions:
\begin{itemize}
\item For $a\geq b+c$\footnote{We will allow a matrix $I_d$ to
occur in the list
even if $d=0$. This just means that the row and column
containing such a
symbol is to be removed. In the first item, this means that if
$a=b+c$ we get
right hand 
side which is a $(b+c)\times(b+c)$ matrix. The sign
$\varepsilon_a$ together
with $a-b-c$ is also removed in this case.}:
$${\mathrm
corank}
(w([a,\varepsilon_a];[b,\varepsilon_b];[c,\varepsilon_c]))=
{\mathrm
 corank} 
(w([a-b-c,\varepsilon_a];[b,-\varepsilon_a\varepsilon_b];[c,
-\varepsilon_a\varepsilon_c]))$$
\item For $c\geq a+b$: $${\mathrm
corank}
(w([a,\varepsilon_a];[b,\varepsilon_b];[c,\varepsilon_c]))=
{\mathrm
 corank} 
(w([a,-\varepsilon_c\varepsilon_a];[b,-\varepsilon_c\varepsilon_b];[c-a-b,
\varepsilon_c]))$$
\item For $b+c\geq a>c$: $${\mathrm
corank}
(w([a,\varepsilon_a];[b,\varepsilon_b];[c,\varepsilon_c]))=
{\mathrm
 corank} 
(w([a,-\varepsilon_b\varepsilon_a];[b-(a-c),\varepsilon_b];[c,
-\varepsilon_b\varepsilon_c]))$$
\item For $a+b\geq c>a$: $${\mathrm
corank}
(w([a,\varepsilon_a];[b,\varepsilon_b];[c,\varepsilon_c]))=
{\mathrm
 corank} 
(w([a,-\varepsilon_b\varepsilon_a];[b-(c-a),\varepsilon_b];[c,
-\varepsilon_b\varepsilon_c]))$$
\end{itemize} 

\medskip

There are many more moves which we shall not pursue here. Also,
there is an
evident tensorial nature to the set-up in the sense that a
common factor of
$a,b$, and $c$ will also turn up as a factor in the resulting corank.

\medskip

Now introduce variables $p=a+b$ and $q=b+c$. Then the four moves
above can be
reformulated, in the same order of appearance, and with the same
stipulations,
as follows:
\begin{itemize}
\item $(p,q,b)\to (p-q,q,b)$. 
\item $(p,q,p)\to (p,q-p,b)$.
\item $(p,q,b)\to (q, 2q-p, b-(p-q))$.
\item $(p,q,b)\to (q, 2q-p, b-(q-p))$.
\end{itemize}

Thus, the moves preserve the lattice generated by $(p,q)$ in
${\mathbb Z}^2$,
and the equivalence class of $b$ modulo the ${\mathbb Z}$
lattice
generated by $p$ and $q$.

\subsection{Special case: $b=1$}\label{subsec7}

We focus on the cases that can be computed without using the
moves changing the
$b$. 
It is straightforward to see that after a certain number of $a$
or $c$ moves
(number 1 and 2 on the list) we get (induction), with the
notation from above,
\begin{eqnarray}{\mathrm corank}(w([a,1];[b,1];[c,1]))=&\\\nonumber
{\mathrm
corank}
(w([X_Lp-Y_Lq-b,(-1)^{X_R+Y_R-1}];[b,(-1)^{X_L+X_R+Y_L+Y_R}];&\\\nonumber
\left[Y_Rq-X_Rp-b,(-1)^{X_L+Y_L-1}\right]))&.\end{eqnarray}

Set $\Delta=g.c.d.(p,q)$. Assume $\Delta>1$. Then we can only
end in a
configuration $(a',1,a')$ with $\Delta=a'+1$. To have a corank
equal to
$a'=\Delta-1$ we need the signs $(-1)^{X_R+Y_R-1}$ and
$(-1)^{X_L+Y_L-1}$ to
equal. This happens exactly when $X_L+X_R+Y_L+Y_R$ is even. If
it is odd, then
we get a $(-1)$ with the middle piece which cancels the $b=1$,
so in this case,
the corank is 1. Write $p=x\Delta$ and $q=y\Delta$. Then the
condition can be
stated as follows:
\begin{eqnarray}
x+y \textrm{ even}: \textrm{corank}&=&\Delta-1,\\\nonumber
x+y \textrm{ odd}: \textrm{corank}&=&1.
\end{eqnarray}

\medskip

Now assume
$\Delta=1$. Then we end up with a configuration
$X_Lp-Y_lq=1=Y_Rq-X_Rp$, though
the two equations, naturally, need not appear simultaneously.
Thus
$X_L+X_R+Y_L+Y_P$ has the same parity as $p+q$ and we get a
non-trivial corank
exactly if $p+q$ is odd. Here, the corank is 1.

Thus, we can state in all cases:

\begin{Prop}The corank $crk(1+\omega^{\mathfrak r})$ of the matrix
\begin{equation}I_{
a+1+c}+\left(\begin{array}{ccc}
0&0&I_c\\0&I_1&0\\I_a&0&0\end{array}
\right)
\end{equation}
is given as follows: Let $p=a+1$, $q=c+1$, let
$\Delta=g.c.d.(p,q)$, and write
$p=x\Delta$ and $q=y\Delta$. Then the condition can be 
stated as follows:
\begin{eqnarray}
x+y \textrm{ even}: crk(1+\omega^{\mathfrak r})&=&\Delta-1,\\
x+y \textrm{ odd}: crk(1+\omega^{\mathfrak r})&=&1.\nonumber
\end{eqnarray}
\end{Prop}

\medskip

\begin{Rem}If we instead start with (for $a>c$) $(a,b+a-c,c)$
then the first $b$
move 
yields $(a,(-1),1,1,c,(-1))$ and all the following $a,c$ moves
preserve these
signs. This more general case is thus covered too.
\end{Rem}

\bigskip

\section{The quantum parabolics}\label{sec8}

\subsection{Quantized minors}\label{subsec8}

Following a construction of classical minors by S. Fomin and
A. Zelevinsky
\cite{fz}, the last mentioned and A. Berenstein have introduced a
family of quantized
minors $\Delta_{u\cdot\lambda,v\cdot\lambda}$ in \cite{bz}. These
are elements of
the quantized coordinate ring ${\mathcal O}_q(G)$. 

{The element 
$\Delta_{u\cdot\lambda,v\cdot\lambda}$ is determined by two
non-zero vectors
$\underline{w}_u,\underline{w}_v$ of weights $u\cdot\lambda$ and $v\cdot\lambda$
respectively, in the finite-dimensional module determined by the
highest weight
vector $\underline{w}_\lambda$. We will always assume that $u\leq v$. Then
$\Delta_{u\cdot\lambda,v\cdot\lambda}(x)=( \underline{w}_u, x
\underline{w}_v)_\lambda$.

\medskip
The construction in (\cite{bz}) is even more sophisticated than what can be 
glimpsed here, since it is given by ``divided powers'' generators, say 
$X^{(k)}$. Since we are working with $q$ generic, all such powers are set equal 
to (a constant times) $X^k$. 
\medskip

\begin{Lem}{\cite{bz}}The
element $\triangle_{u\lambda,v\lambda}$ indeed depends only on
the weights
$u\lambda,v\lambda$, not on the choices of
$u, v$ and their reduced words.
\end{Lem}

\medskip

\begin{Thm}[A version of Theorem~10.2 in \cite{bz}]
\label{10.2}For any
$\lambda,\mu\in P^+$, and $s, s', t, t' \in W$ such that
$$\ell(s's) = \ell(s') + \ell(s), \ell(t't) = \ell(t') + \ell(t)
,$$the
following holds:
$$ \triangle_{s's\lambda,t'\lambda} \triangle_{s'\mu,t't\mu} =q^{(s\lambda | 
\mu) - (\lambda |
t\mu)}\triangle_{s'\mu,t't\mu} 
\triangle_{s's\lambda,t'\lambda}.$$
\end{Thm}

\bigskip

\bigskip

\begin{Def}
A Levendorskii-Soibelman quadratic algebra is a quadratic algebra with
generators $W_\beta$ labeled by a subset $\Phi^+_0\subseteq
\Phi^+$ of the positive roots together with a total ordering $<$
of $\Phi^+_0$ such that for all
$\beta_1,\beta_2\in\Phi^+_0$,
\begin{equation}\label{7.3}
\beta_1<\beta_2\Rightarrow
W_{\beta_1}W_{\beta_2}=q^{\alpha_{ij}}
W_{\beta_2}W_{\beta_1}\ + \ L.O.T.s.
\end{equation}
The L.O.T.s stand for certain sums of products of 
elements $W_\beta$ where $\beta_1<\beta<\beta_2$, and thus are of lower order 
(hence the abbreviation) in the lexicographical order. The exponents 
$\alpha_{ij}$ are assumed to be integers.
\end{Def}

\begin{Def}Let ${\mathcal A}_q({\mathfrak n}_\omega)\subset {\mathcal 
U}^{\geq}({\mathfrak n}_\omega):={\mathcal U}({\mathfrak n}_\omega)\cdot 
{\mathcal U}^0$ be a
Levendorskii-Soibelman quadratic algebra with a set of generators
$W_\beta$ where $\beta$ runs through the roots of ${\mathfrak
u}={\mathfrak n}_\omega$ and ordered according to $\omega$. We will say that 
${\mathcal A}_q({\mathfrak n}_\omega)$ splits ${\mathcal
U}^{\geq}({\mathfrak n}_\omega)$, or that we have a splitting furnished by 
${\mathcal
A}_q({\mathfrak n}_\omega)$ if
\begin{equation}
 {\mathcal U}^{\geq}({\mathfrak n}_\omega)={\mathcal A}_q({\mathfrak 
n}_\omega)\times_s {\mathcal U}^{0}.
\end{equation}
\end{Def}

\medskip

We now return to the setting of Section~\ref{sec5} and consider an element 
$\omega^{\mathfrak r}\in W^{\mathfrak p}$.

\begin{Def}${\mathcal M}_q({\mathfrak r})$ denotes the algebra generated by 
the previously defined elements $Z_{\beta_j}\equiv Z_{s,t}$.\end{Def}

It is well known that this is a 
Levendorskii-Soibelman quadratic algebra (\cite{lev}). More detailed proofs 
have been given in \cite{cp} and \cite{jan2}. This is equal to the algebra in 
\cite{leclerc} except $q\to q^{-1}$.

One considers in \cite{leclerc}, and transformed to our terminology, more 
general elements 
\begin{equation}\label{59}D_{\xi,\eta}=\triangle_{\xi,\eta}K^{-\eta}.\end{equation}

The family $D_{\xi,\eta}$
satisfies equations analogous to those in  Theorem~\ref{10.2} subject to the 
same restrictions on the relations between the weights. 

The elements 
\begin{equation}W_{s,t}=Z_{s,t}K^{\omega_{s,t}(\Lambda_s)}\end{equation} thus 
generate a 
quadratic algebra.

\begin{Rem}
More generally, one may consider algebras generated by elements $$ X^{\overline F}_\gamma=Z_\gamma K^{\overline{F}(\gamma)}$$where $\overline{F}(\gamma)$ is in the root lattice. Not all such algebras are of interest. A natural requirement to have fulfilled by $\overline F$ is that the elements $ X^{\overline F}_\gamma$ satisfy equations analogous to (\ref{7.3}). An analogous approach is to consider elements $\Delta^F_{\xi,\eta}= \Delta_{\xi,\eta} K^{F(\xi,\eta)}$ where, for all $u,v\in W$ and $i\in I$ with $l(us_i)=l(u)+1$ and
$l(vs_i)=l(v)+1$,
\begin{eqnarray}\label{66}
F(us_i(\Lambda_i),vs_i(\Lambda_i))+F(u(\Lambda_i),v(\Lambda_i))&=&\\\nonumber
F(us_i(\Lambda_i),v(\Lambda_i))+F(u(\Lambda_i),
vs_i(\Lambda_i))&=&\\
\sum_{j\neq i}-a_{ji}F(u(\Lambda_j),v(\Lambda_j)).\nonumber
\end{eqnarray}
For fixed integers $n,m$ the function$$F^{n,m}_{\xi,\eta}=n\xi+m\eta$$meets the requirements in (\ref{66}).

Similarly, any linear function $\overline F$ yields a Levendorskii-Soibelman  quadratic algebra as follows by considering weights in (\ref{7.3}). 

Some algebraic properties of different such algebras, eg. their P.I. degree, may be different, see \cite{jz2}. 
\end{Rem}

\begin{Def}${\mathcal W}_q({\mathfrak r})$ denotes the algebra generated by the 
elements $W_{s,t}$.\end{Def} 

Its center will be discussed later.

\medskip

Furthermore, a quasi-polynomial algebra $\overline{\mathcal W}_q({\mathfrak 
r})$ is associated to  ${\mathcal W}_q({\mathfrak r})$ in the following obvious 
way: The generators are
\begin{equation}w_{s,t}=z_{s,t}K^{\omega_{s,t}(\Lambda_s)},\end{equation} where 
the elements $z_{s,t}$ 
are the old ones of Section~\ref{sec6}.

The relations are easily 
\begin{equation}
(s,t)<(c,d)\Rightarrow w_{s,t}w_{c,d}=q^{R_{s,t}^{c,d}}w_{c,d}w_{s,t},
\end{equation}
where
\begin{equation}R_{s,t}^{c,d}=( \beta_{s,t},\beta_{c,d})+(
\Lambda_s,\beta_{c,d})-(\Lambda_c,\beta_{s,t})-(\sum_{i=1}^t\beta_{s,i},\beta_{c
,d})+(\sum_{j=1}^d\beta_{c,j},\beta_{s,t}).
\end{equation}

\bigskip

We end this section with the following easy consequence of Theorem~\ref{10.2}:

\begin{Prop}Each element $\triangle_{\Lambda_s,\omega^{\mathfrak r}(\Lambda_s)}, 
s\in Im(\pi_{\mathfrak r})$, quasi-commutes with all elements of the form 
\begin{equation}\triangle_{\omega_{c,d}(\Lambda_c),\omega_{c,f}(\Lambda_c)}.
\end{equation}
A similar result holds for the elements $D_{\Lambda_s,\omega^{\mathfrak 
r}(\Lambda_s)}$.
\end{Prop}

It is proved in (\cite{leclerc}) that the elements $Z_{a,b}$ are of the form 
\begin{equation}D_{\omega_{c,d}(\Lambda_c),\omega_{c,f}(\Lambda_c)}
\end{equation}
for some elements $(c,d),(c,f)\in{\mathbb U}({\mathfrak r})$, but we actually only need that the elements $Z_{a,b}$ are rational functions in the elements $D_{\Lambda_s,\omega_{c,s}(\Lambda_s)}$. This holds in general.  Hence we conclude 
in particular, cf. (\ref{59}), that

\begin{Prop}\label{quasipol}The elements $C_{s}({\mathfrak 
r})=D_{\Lambda_s,\omega^{\mathfrak r}\Lambda_s}$ $q$-commute with all elements 
in the algebra ${\mathcal M}_q({\mathfrak n}_\omega)$.
\end{Prop}

\medskip
\begin{Def}An element $C\in{\mathcal M}_q({\mathfrak n}_\omega)$ that $q$-commutes with all elements 
in the algebra ${\mathcal M}_q({\mathfrak n}_\omega)$ is said to be covariant.
\end{Def}

\bigskip

\section{The centers of the nilpotent part of quantized parabolics}\label{sec9}

It follows easily from Theorem~\ref{10.2}  that we have the following formula:

\begin{Prop}[\cite{bz},\cite{leclerc}]Let $(s,t)\leq (c,d)$, then
\begin{equation}\label{lec}D_{\Lambda_s,\omega_{s,t}\Lambda_s}D_{\Lambda_c,
\omega_{c,d}\Lambda_c}=q^{((1-\omega_{s,t})(\Lambda_{s}),(1+\omega_{c,d}
)(\Lambda_{c}))}D_{\Lambda_c,
\omega_{c,d}\Lambda_c}D_{\Lambda_s,\omega_{s,t}\Lambda_s}.\end{equation}
\end{Prop}
There is a sign difference in the $q$ exponent in relation to 
Corollary~\ref{only}, otherwise they are identical.

\medskip

Secondly, the critical step in the proof of Proposition~\ref{key} is where one 
considers $\omega_{s,t}=\omega'\circ\omega_s$ in the pairs \begin{equation}
(\Lambda_s,\omega'\Lambda_s),((\Lambda_s,\omega'\circ\sigma_s\Lambda_s),a_{sj}
(\Lambda_j,\omega'\Lambda_j) \ (\textrm{for all }a_{sj}<0).
\end{equation}
Combined with equation (\ref{5.11-f}) we obtained the important 
Proposition~\ref{6.16}. Now we have equation (\ref{lec}) which is the same hence 
furnishes  the same conclusions.

\medskip

Thirdly, the elements $D_{\Lambda_s,\omega_{s,t}\Lambda_s}$ are ordered 
in the same way as the elements $\Delta_{\Lambda_s,\omega_{s,t}\Lambda_s}$, and 
as 
the elements $Z_{s,t}$ and $z_{s,t}$.

\medskip
Fourthly, it follows easily from \cite[Lemma~11.4 and Corollary~12.4]{leclerc} 
that the elements $C_{s}({\mathfrak r})$ are regular and 
that all elements $M^\downarrow_{c,d}=D_{\Lambda_c,\omega_{c,d}\Lambda_c}$ are 
polynomials in the generators $Z_{s,t}$. We also know by 
proposition~\ref{quasipol} that they generate a quasi-polynomial subalgebra. We 
can then invert them and consider now the algebra
\begin{equation}
\tilde{\mathcal M}_q({\mathfrak n}_\omega)={\mathcal M}_q({\mathfrak 
n}_\omega)[C_{1}({\mathfrak r})^{-1},\dots, C_{R}({\mathfrak r})^{-1}].
\end{equation}

\medskip

We can then state, cf. Remark~\ref{ofcourse},

\begin{Thm}[\cite{yak}]The center of the algebra $\tilde{\mathcal M}_q({\mathfrak 
n}_\omega)[C_{1}({\mathfrak r})^{-1},\dots, C_{R}({\mathfrak r})^{-1}]$ is 
given 
by the kernel of 
$(1+\omega^{\mathfrak r})$.
\end{Thm}

\bigskip

\subsection{Example: The full nilpotent - the longest element $\omega_0\in 
W$.}\label{subsec9}

We can now revisit Subsection~\ref{subsec62}, using the same numbering, and easily obtain:

\begin{Prop}[\cite{Ca0}]
If the simple Lie algebra is not in the list (\ref{exc1}), the center of 
${\mathcal U}_q({\mathfrak n})$ is 
generated by the elements $C_{s}({\mathfrak r})$, with $s=1\dots,R$.
In type ${\mathbf A}_\ell$, with simple roots
$\alpha_1,\alpha_2,\dots,
\alpha_\ell$ we have that 
\begin{equation}\omega_m(\Lambda_{\alpha_s})=-\Lambda_{\alpha_{\ell-s+1}},\end{equation}
and the center is generated by the elements\begin{eqnarray*}& 
C_{s}({\mathfrak r})C_{\ell-s+1}({\mathfrak r})\textrm{ for 
}2s\neq\ell+1\\
&C_{s}({\mathfrak r})\textrm{ for }2s=\ell+1\end{eqnarray*}
In type ${\mathbf D}_{2\ell+1}$, the center is generated by the 
elements\begin{equation}C_{s}({\mathfrak 
r})(s=1,\dots,2\ell-1),\textrm{ and }C_{2\ell}({\mathfrak r}) 
C_{2\ell+1}({\mathfrak r}).\end{equation}
In type ${\mathbf E}_6$, the center is generated by the elements 
\begin{equation}
C_{1}({\mathfrak r})C_{6}({\mathfrak r}),  
C_{2}({\mathfrak r})C_{5}({\mathfrak r}), C_{2}({\mathfrak 
r}),\textrm{ and }C_{4}({\mathfrak r}).\end{equation}
\end{Prop}

\bigskip

\begin{Ex}For the $a\times b$ quantized matrix algebra matrix we get, provided 
$a>1$ and $b>1$,
\begin{equation}\dim({\mathcal C}_{a\times b})=2 
+\left(g.c.d.(a-1,b-1)-1\right).
\end{equation}
In case $a\geq b=1$, the algebra is commutative; $\dim({\mathcal 
C}_{a\times1})=a$.
\end{Ex}

\bigskip

We now turn our attention to the Berenstein-Zelevinski algebra ${\mathcal 
W}_q({\mathfrak r})$.

We  will for simplicity assume that
\begin{equation}\label{66}
\omega^{\mathfrak r}=\omega_1\omega_2\dots\omega_\ell,
\end{equation}
where
\begin{equation}
\forall i,j :\ Z_{i,j}=-\omega_1\dots\omega_j(\alpha_i).
\end{equation}

This is satisfied for all types except $D_n, E_6,E_7,E_8$.

\smallskip

The precise meaning of this condition is that each $\omega_i$ can be written as 
a product of pairwise different reflexions $\sigma_{i,k}$. Furthermore, if 
$\sigma_{i,k}$ takes part in $\omega_i$ then it also takes part in 
$\omega_{i-1}$ (etc.).  So, in particular, $Im(\pi_{\omega_i})\subseteq 
Im(\pi_{\omega_{i-1}})$.

This  gives a natural ordering, for each $j$, of the elements $Z_{i,j}$ as well 
as the elements 
\begin{equation}
 \triangle_{s,t}=\Delta_{\Lambda_s,\omega_{s,t}\Lambda_s},
\end{equation}
where $\omega_{s,t}=\omega_1\dots\omega_t$. We also 
remark that if $\omega_1$ has $\sigma_x$ to the right; 
$\omega_1=\omega^{(1)}\sigma_x$ with $\omega^{(1)}$ reduced, then 
$\Delta_{\Lambda_x,\omega_1\Lambda_x}$ commutes with all other elements 
$\triangle_{s,t}$. Any rewriting 
of $\omega_1$ that results in another reflection $\sigma_y$ at the right end 
will similarly give rise to a central element.

\medskip

We  see immediately that we have the following result:

\begin{Prop}
Under the above assumption on $\omega^{\mathfrak r}$, there is  a non-trivial 
center.
\end{Prop}

For the algebras $\tilde{\mathcal M}_q({\mathfrak n}_\omega)$ it may happen that 
there is a trivial center. 

\medskip

Let ${\mathcal L}_\triangle({\mathfrak r})$ denote the Laurent quasi-polynomial 
algebra generated by the elements in the set $\{\triangle_{s,t}\mid (s,t)\in 
{\mathbb U}({\mathfrak r})\}$, and let ${\mathbb L}_\triangle$ denote the 
symplectic form defined by these elements. We may also assume that there is an 
auxiliary inner product in which these elements form an orthonormal set. 

A very important difference between the current case and the previous is that 
the elements $\triangle_{\Lambda_s,\Lambda_s}$ are not equal to the constant 
function.
This gives some modifications. However, using the assumption (\ref{66}) it is 
straightforward to see that one may make a construction similar to the case of 
$\overline{\mathcal L}_q({\mathfrak r})$, especially Definition~\ref{6.12}, to 
obtain:

\medskip

\begin{Lem}\label{9.6}Let $(i,j_0)$ be given with $j_0\neq1$. Then there exists 
an element $\nabla_{i,j_0}$ in ${\mathcal L}_\triangle({\mathfrak r})$ such that 
\begin{equation}\triangle_{s,t}\nabla_{i,j_0}=q^{\Gamma_{i,j_0,s,t}}\nabla_{i,
j_0}    \triangle_{s,t}
\end{equation}with 
\begin{equation}
\Gamma_{i,j_0,s,t}=\left\{\begin{array}{rl}0&\textrm{ if }s\neq i\\-{d_i}&\textrm{ 
if }s=i\textrm{ and }j_0\leq t\leq 
i_{max}\\+{d_i}&\textrm{ if }s=i\textrm{ and }t<j_0\end{array}\right. .
\end{equation}
\end{Lem}

By considering elements of the form $\nabla_{s,t}\nabla_{s,t-1}^{-1}$  we get ($q$ generic)

\begin{Lem}\label{9.7}
\begin{equation}\forall i,\forall 1<j_0<i_{max}: \triangle_{i,j_0}\in 
\textrm{Range}({\mathbb L}).
\end{equation}
\end{Lem}

Combining Lemmas \ref{9.6} and \ref{9.7} we also get

\begin{Lem}\label{9.8}
\begin{equation}\forall i: \triangle_{i,1}^{-1}\triangle_{i,i_{max}}\in 
\textrm{Range}({\mathbb L}).
\end{equation}
\end{Lem}

This, similarly to the previous case, easily implies

\begin{Lem}\label{9.9} The center is generated by elements of the form 
\begin{equation}\prod_i(\triangle_{i,1}\triangle_{i,i_{max}})^{n_i}.
\end{equation}
\end{Lem}

\medskip

It is understood that in case $\triangle_{i,1}=\triangle_{i,i_{max}}$ we omit 
one of these factors. 

\smallskip

Naturally, a central element must commute with all elements $\triangle_{s,1}$. 
On the other hand, it follows that it commutes with all elements of the form 
$\nabla_{s,t_0}$ ($t_0>1$).
Thus, it is also sufficient that it commutes with the elements 
$\triangle_{s,1}$. In other words, the following must hold:

\begin{Lem}The element $\prod_i(\triangle_{i,1}\triangle_{i,i_{max}})^{n_i}$ is 
central if and only if
\begin{equation}
 \forall s: \left(\Lambda_s,(1-\omega_1^{-1}\omega)\sum_i n_i\Lambda_i\right)=0.
\end{equation}
\end{Lem}

Consider the algebra
\begin{equation}
\tilde{\mathcal W}_q({\mathfrak n}_\omega)={\mathcal W}_q({\mathfrak 
n}_\omega)[(\triangle_{1,1}\triangle_{1,1_{max}})^{-1},\dots, 
(\triangle_{R,1}\triangle_{R,R_{max}})^{-1}].
\end{equation}

\medskip

We can then state

\begin{Thm}The center of the algebra $\tilde{\mathcal W}_q({\mathfrak r})$ is 
given 
by the kernel of 
$(1-\omega_1^{-1}\omega^{\mathfrak r})$ on $S_{\mathfrak r}$.
\end{Thm}

\bigskip

Let $\check C_{i}({\mathfrak r}), i=1,\dots,R$ denote the elements in ${\mathcal 
U}_q({\mathfrak n}_\omega)$ that correspond to the functions $C_{i}({\mathfrak 
r}), i=1,\dots,R$.

In view of our discussion in Section~\ref{sec2} we can also state for any 
$\omega\in W$:

\begin{Thm}
The center of the algebra ${\mathcal U}_q({\mathfrak n}_\omega)[\check 
C_{1}({\mathfrak r})^{-1},\dots, \check C_{R}({\mathfrak r})^{-1}]$ is given by
by the kernel of 
$(1+\omega)$.
\end{Thm}

\bigskip

For our final class of examples, let $\omega^{\mathfrak a},\omega^{\mathfrak c}$ be elements in $W^{\mathfrak p}$ such that $\omega^{\mathfrak a}<\omega^{\mathfrak c}$.  Obviously, 
the elements in 
\begin{equation}
{\mathcal M}_q(\omega^{\mathfrak a},\omega^{\mathfrak c})=\{Z_{s,t}\mid s_{\mathfrak a}\leq t\leq s_{\mathfrak c}\}
\end{equation} are the generators of a quadratic algebra which is a subalgebra of  ${\mathcal M}_q(\omega^{\mathfrak c})$.
Hence the algebra ${\mathcal A}(\omega^{\mathfrak a},\omega^{\mathfrak c})$ generated by these elements is a Levendorskii-Soibelman quadratic algebra. We may extend the definition of the previous elements $C_s({\mathfrak r})$ to 
\begin{equation}C_s({\mathfrak a},{\mathfrak c})=D_{\omega^{\mathfrak a}\Lambda_s,\omega^{\mathfrak c}\Lambda_s}
\end{equation} 
and similarly introduce an algebra
\begin{equation}
\tilde{\mathcal M}_q({\mathfrak a},{\mathfrak c})={\mathcal M}_q({\mathfrak 
a},{\mathfrak c})[C_{1}({\mathfrak a},{\mathfrak c})^{-1},\dots, C_{R}({\mathfrak a},{\mathfrak c})^{-1}].
\end{equation}

\begin{Prop}The center  of $\tilde{\mathcal M}_q({\mathfrak a},{\mathfrak c})$ is given by the null space of $\omega^{\mathfrak a}+\omega^{\mathfrak c}$.
\end{Prop}

\bigskip

\section{Root of unity}

We here briefly discuss how our constructions and results may be specialized to an $\varepsilon\in{\mathbb C}$ for which $\varepsilon^m=1$ for some positive integer $m$. We define our quantized enveloping algebra ${\mathcal U_\varepsilon}({\mathfrak g})$ as in \cite{cp}[\S19]. 
\medskip

First of all then, the associated quasi-polynomial algebra is taken to be 

$\overline{\mathcal U}_\varepsilon({\mathfrak n}_\omega)$ with 
generators $z_1,\dots,z_{\ell(w)}$, and  relations 
\begin{equation}z_{i}z_{j}=\varepsilon^{(\gamma_i,\gamma_j)}z_{j}z_{i}.\end{equation}

\smallskip

Secondly, the result of Levendorskii and Soibelman (\cite{lev2}, \cite{lev0})  hold also in this generality so that we have a quadratic algebra ${\mathcal M}_\varepsilon(\omega^{\mathfrak r})$ generated by elements $Z_1,\dots,Z_{\ell(w)}$, and  relations 
\begin{equation}Z_{i}Z_{j}=\varepsilon^{(\gamma_i,\gamma_j)}Z_{j}Z_{i} + L.O.T..\end{equation}

Our minors can be (and have been) rescaled so that they are polynomials in these generators. Thus we get the same by directly specializing as we get from specializing the Lusztig form.

We  assume that $m$ is odd and, in case there is an irreducible component  of $G_2$, $m$ should be relatively prime to $3$. If $k$ is any field of characteristic ${\it char}(k)\neq 2,3$ our results even hold here.

First of all then, all elements $z^m_i,Z^m_i$ are central in the appropriate algebras (\cite{lus-c1}, \cite{lus-fin}).

More generally, the following holds:

\begin{Prop}
Any $D^m_{\omega_{s,d}\Lambda_s,\omega_{s,t}\Lambda_s}$ ($d<t\leq s_{\mathfrak r}$)  is central. In particular it holds for any $s$ that  $C_s({\mathfrak r})^m$ is central.
\end{Prop}
This follows from (\ref{5.11}) in combination with (\cite{leclerc}).

After that we can make replacements \begin{equation}C_i({\mathfrak r})^{-1}\rightarrow C_i({\mathfrak r})^{m-1}\end{equation}for any $i$. 
In particular, there is no need to localize, and the above mentioned elements are central in ${\mathcal M}_\varepsilon(\omega^{\mathfrak r})$.

With these remarks in mind, we need then only focus on the block diagonal form of ${\mathcal L}$. We take over terminology from (\cite{jaz1}, \cite{jj}, \cite{jp}).  The main task is to determine the block diagonal form of ${\mathcal L}_\varepsilon$. In the mentioned references it is also seen that the finer details of such block diagonal forms can be quite delicate.  Recall

\begin{Prop}[\cite{cp}] Let $h$ denote the cardinality of the map  (we maintain the notation)  
\begin{equation}
{\mathcal L}_\varepsilon: {\mathbb Z}^r\rightarrow ({\mathbb Z}/m{\mathbb Z})^r.\end{equation}Then the P.I. degree of ${\mathcal M}_\varepsilon(\omega^{\mathfrak r})$ is $\sqrt{h}$.
\end{Prop}

We wish to discuss this cardinality.

\medskip 
For this purpos observe that (the analogue of) Proposition~\ref{6.16} implies that one can bring the matrix  ${\mathcal L}$ into a column reduced form 
\begin{equation}\label{93}
\left(\begin{array}{cc}2I_{r-r_0}& Y\\0&Z\end{array}\right),
\end{equation}
where $r_0=\vert Im(\pi_{\mathfrak r})\vert$.

The matrix $Z$ has entries of the form $\left((1+\omega^{\mathfrak r})\Lambda_s, (1-\omega_{c,d})\Lambda_c\right)$ as $c,s\in Im(\pi_{\mathfrak r})$. Typically, $d=c_{\mathfrak r}$ or $d=c_{{\mathfrak r}-1}$, but we refrain from pursuing all details here.

We immediately get

\medskip

\begin{Prop}
The $q$-generic center is trivial if and only if $\det(Z)\neq0$. In this case, $\det({\mathcal L})=2^{r-r_0}\det(Z)$. In particular, if  
\begin{equation}
\left(\begin{array}{cc}0& p\\-p&0\end{array}\right)
\end{equation}
is a member of a block diagonal form of ${\mathcal L}$ then $p^2$ divides $2^{r-r_0}\det(Z)$. Furthermore, if $p_1^t$ is a divisor in $\det(Z)$, with $p_1\neq2$ a prime, and $t$ maximal, then $t$ must be even and the block diagonal contains factors
\begin{equation}
\left(\begin{array}{cc}0& u_i p_1^{s_i}\\-u_i p_1^{s_i}&0\end{array}\right)
\end{equation}
such that $2\sum_is_i=t$, and $u_i,s_i\in{\mathbb N}$ with each $u_i$ relatively prime to $p_1$.
\end{Prop}

Further observations about the elements in the block diagonal can be deduced. However, in case the center of $(1+\omega^{\mathfrak r})$ is nontrivial, it is more complicated to make substantial claims without invoking specific properties of $\omega^{\mathfrak r}$: Let ${\mathbf d}=R-i_0$ be the dimension of said null space. Consider the change of basis map $T$, \begin{eqnarray}&T: \{C_1({\mathfrak r}), \dots, C_{i_0}({\mathfrak r}),C_{i_0+1}({\mathfrak r}),\dots, C_r({\mathfrak r})\}\rightarrow\\\nonumber   &\{C_1({\mathfrak r}), \dots, C_{i_0}({\mathfrak r}),U_{1}({\mathfrak r})\dots, U_{\mathbf d}({\mathfrak r})\}
\end{eqnarray}
where $\{U_{1}({\mathfrak r})\dots, U_{r-i_0}({\mathfrak r})\}$ is a basis of the center. Assume that $T$ as a linear map has an integer matrix of determinant $1$. We then have

\begin{equation}
{\mathcal L}=A^t\left(\begin{array}{cc}BD&0 \\0&0\end{array}\right)A
\end{equation}
where $BD$ is a block diagonal $(r-{\mathbf d})\times (r-{\mathbf d})$ matrix in which each block is a skew symmetric $2\times 2$ matrix, and where 
\begin{equation}
A=\left(\begin{array}{cc}A_1& 0\\0&1_{\mathbf d}\end{array}\right)
\end{equation}with $\det(A_1)=1$.

Proposition~\ref{6.16}  can then be used again, but we need to use the matrix $T^{-1}$ to replace factors of $C_i({\mathfrak r})$, $i=i_0+1,\dots, r$ in ${\mathbb B}$ by extra factors (integer powers) coming from the set $\{C_1({\mathfrak r}), \dots, C_{i_0}({\mathfrak r})\}$. After that, one may again  perform column operations by integer matrices, but there may be determinants other than 1.

It follows that 
\begin{equation}
\det(BD)\det(A_2)=2^{r-{\mathbf d}} \det(Z_1)
\end{equation}
where $A_2$ is an invertible integer matrix whose determinant depends on the finer details of the change-of-basis.

Furthermore, $Z_1$ may be construed as the previous $Z$ modulo the null space of $1+\omega^{\mathfrak r}$.

Thus we obtain

\begin{Prop}Under the above hypothesis, 
if a  prime $p>2$ divides $\det(BD)$ then it also divides  $\det(Z_1)$. 
\end{Prop}

\end{document}